\documentclass{amsart}
\usepackage{amssymb}
\usepackage{graphicx}
\usepackage{url}
 \newtheorem{theorem}{Theorem}
 \newtheorem{lemma}[theorem]{Lemma}
 \newtheorem{proposition}[theorem]{Proposition}
 \newtheorem{corollary}[theorem]{Corollary}
 \newtheorem{remark}[theorem]{Remark}
 \newtheorem{example}[theorem]{Example}
 \newtheorem{definition}[theorem]{Definition}
 \newtheorem{conjecture}[theorem]{Conjecture}

 \newtheorem{question}[theorem]{Question}

\newcommand{\bpr}{\begin{proof}}
\newcommand{\epr}{\end{proof}}
\newcommand{\beq}{\begin{equation}}
\newcommand{\eeq}{\end{equation}}
\newcommand{\bThm}{\begin{theorem}}
\newcommand{\eThm}{\end{theorem}}
\newcommand{\blem}{\begin{lemma}}
\newcommand{\elem}{\end{lemma}}
\newcommand{\bpro}{\begin{proposition}}
\newcommand{\epro}{\end{proposition}}
\newcommand{\bcor}{\begin{corollary}}
\newcommand{\ecor}{\end{corollary}}
\newcommand{\brem}{\begin{remark}}
\newcommand{\erem}{\end{remark}}
\newcommand{\bexa}{\begin{example}}
\newcommand{\eexa}{\end{example}}
\newcommand{\bdf}{\begin{definition}}
\newcommand{\edf}{\end{definition}}
\newcommand{\bcon}{\begin{conjecture}}
\newcommand{\econ}{\end{conjecture}}
\newcommand{\bque}{\begin{question}}
\newcommand{\eque}{\end{question}}

 \newcommand{\R}{{\mathbb R}}

\newcommand{\comment}[1]{}

\title{Some generalizations of Satoh's Tube map}
\author{Blake K. Winter}
\address{}

\date{}

\begin{document}
\thispagestyle{empty}

\begin{abstract}
Satoh has defined a map from virtual knots to ribbon surfaces embedded in $S^4$. Herein, we generalize this map to virtual $m$-links, and use this to construct generalizations of welded and extended welded knots to higher dimensions. This also allows us to construct two new geometric pictures of virtual $m$-links, including 1-links.
\end{abstract}

\maketitle

%
\section{Introduction}
%

In \cite{SS}, Satoh defined a map from virtual link diagrams to ribbon surfaces knotted in $S^4$ or $\R^4$. This map, which he called $Tube$, preserves the knot group and the knot quandle (although not the longitude). Furthermore, Satoh showed that if $K\cong K'$ as virtual knots, $Tube(K)$ is isotopic to $Tube(K')$. In fact, we can make a stronger statement: if $K$ and $K'$ are equivalent under what Satoh called $w$-equivalence, then $Tube(K)$ is isotopic to $Tube(K')$ via a stable equivalence of ribbon links without band slides. For virtual links, $w$-equivalence is the same as welded equivalence, explained below. In addition to virtual links, Satoh also defined virtual arc diagrams and defined the $Tube$ map on such diagrams.
On the other hand, in \cite{BKW, BKW2}, the notion of virtual knots was extended to arbitrary dimensions, building on the geometric description of virtual knots given by Kuperberg, \cite{Kup}. It is therefore natural to ask whether Satoh's idea for the $Tube$ map can be generalized into higher dimensions, as well as inquiring into a generalization of the virtual arc diagrams to higher dimensions.
In this paper, we will accomplish both of these things. First, we will review the definition of virtual $n$-knots and Satoh's $Tube$ map for virtual 1-links and arcs. Then, we will generalize virtual $m$-knots to include a type of "virtual $m$-knotoid." We will then define a map $T_{mn}$ on virtual $m$-links/arcs. This map will define an $m$-link, and the isotopy class of the resulting $m$-link is independent of the choices involved as long as $m\geq 2n$. The $T_{mn}$ map also preserves the fundamental group and quandle. As a bonus, the definition of this map leads to a natural generalization of Satoh's $w$-equivalence for $m$-knots. Since this is the same as welded equivalence for virtual $1$-links, this provides a natural generalization of welded equivalence into higher dimensions. This also gives two other geometric interpretation of virtual $n$-links, including 1-links (that is, even for 1-links, this seems to be a somewhat new interpretation of them).

\section{Virtual links and \emph{w}-equivalence}
The theory of virtual links was developed by Kauffman, \cite{Kauff}, as a generalization of classical links. His original definition was combinatorial in nature, but we will begin with the geometric definition, due to Kuperberg, \cite{Kup}, and others, \cite{KamaKama, StableEq}.

Every classical $1$-link can be expressed by a \emph{link diagram}, which is a combinatorial object consisting of a planar graph on a surface whose vertices are $4$-valent. At each vertex, two of the edges of the graph are marked as \emph{the overcrossing arc} and the other two are marked as the \emph{undercrossing}. These markings give information on how to resolve them in three dimensions.

Let  $K:M\hookrightarrow F\times [0,1]$ be a $1$-link, where $F$ is a surface (possibly with boundary) and $M$ is a disjoint union of circles. Then there are canonical projections $\pi:F\times [0,1]\rightarrow F$ and $\pi_I:F\times [0,1]\rightarrow [0,1]$. It is possible to isotope $K:M\hookrightarrow F\times I$ such that $\pi \imath$ is an immersion and an embedding except for a finite number of double points called \emph{crossings}. Note that we will often abuse notation by using $K$ to refer to both the embedding and its image.
In order to recover $K$, at each double point, we mark which strand projects under $\pi_I$ to the smaller coordinate in $[0,1]$ at that double point; this is the undercrossing. These marks allow us to recover $K$ up to ambient isotopy. Observe that $\pi(K)$ is a $4$-valent graph. The interiors of the edges of this graph will be called \emph{semi-arcs} of the link diagram. Let $D_{-}$ be the set $\{x\in K | \pi^{-1}(\pi(x))=\{x, y\}, x \neq y, \pi_I(x)<\pi_I(y)\}$. For any connected component $A$ of $K-D_-$, $\pi(A)$ is called an \emph{arc} of the link diagram.

For diagrams of $1$-links, there exists a set of three moves, the \emph{Reidemeister moves}, which satisfy the following condition: if $K$ and $K'$ are ambient isotopic links in $F\times [0,1]$, then the diagrams of $K$ and $K'$ on $F$ differ by a finite sequence of ambient isotopies of $F$ and Reidemeister moves. Thus, questions about $1$-links may be translated into questions about equivalence classes of link diagrams, with the equivalence relation being generated by the Reidemeister moves.

We now wish to define \emph{virtual equivalence}. Define a relation $\sim$ by the condition that a link $K_1$ in $F_1 \times I$ is $\sim$-related to $K_2$ in $F_2\times I$ iff there exists an embedding $f:F_1\rightarrow F_2$ such that $f\times id_{[0,1]}(K_1)=K_2$. We define virtual equivalence to be the equivalence relation generated by $\sim$ together with smooth isotopy links in thickened surfaces.

We now wish to relate this to the combinatorial definition of Kauffman. Since a link diagram is inherently a graph on a surface $F$, each link diagram gives rise to a \emph{Gauss code}. There are various ways of defining Gauss codes for links. All of them specify the set of arcs in the diagram, and, for each crossing, the two arcs which terminate at that crossing are specified, together with the order in which these crossings occur on their overcrossing arcs. Concretely, a Gauss diagram for $K:S^1\cup ...\cup S^1\xrightarrow F\times I \xrightarrow{\pi} F$ is $K$ together with a collection of triples $(x,y,z)$, one for each vertex $v$ of the diagram $\pi\imath(K),$ where $x,y$ are preimages of $v$ under $\pi$, $x$ is the overcrossing (i.e. $\pi_I(y)>\pi_I(x)$), and $z=\pm 1$, depending on whether the crossing is right or left-handed. These triples are usually indicated by arrows whose tails are at the overcrossing point, and whose heads are at the undercrossing point. The right-handed crossings may be labeled with a + and the left-handed crossings with a -, rather than using $\pm 1$. 

\begin{figure}
		\centering
			\includegraphics[scale=0.3]{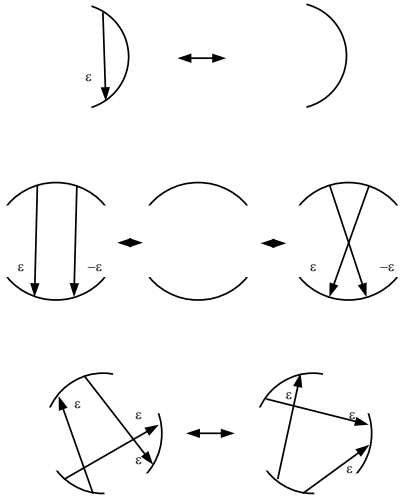}
		\caption{The effects of the three Reidemeister moves on a Gauss code. Here $\epsilon$ denotes the sign of the crossing and $-\epsilon$ denotes a crossing of the opposite sign.}
		\label{gaussreid}
\end{figure}

However, not every Gauss code corresponds to a link diagram on $\R^2$. A Gauss code which corresponds to a link diagram on $S^2$ is termed \emph{realizable}. One method for defining virtual links is to consider arbitrary Gauss codes, modulo local changes corresponding to Reidemeister moves. These moves for Gauss codes are illustrated in Fig. \ref{gaussreid}. 

\begin{figure}
		\centering
			\includegraphics[scale=0.5]{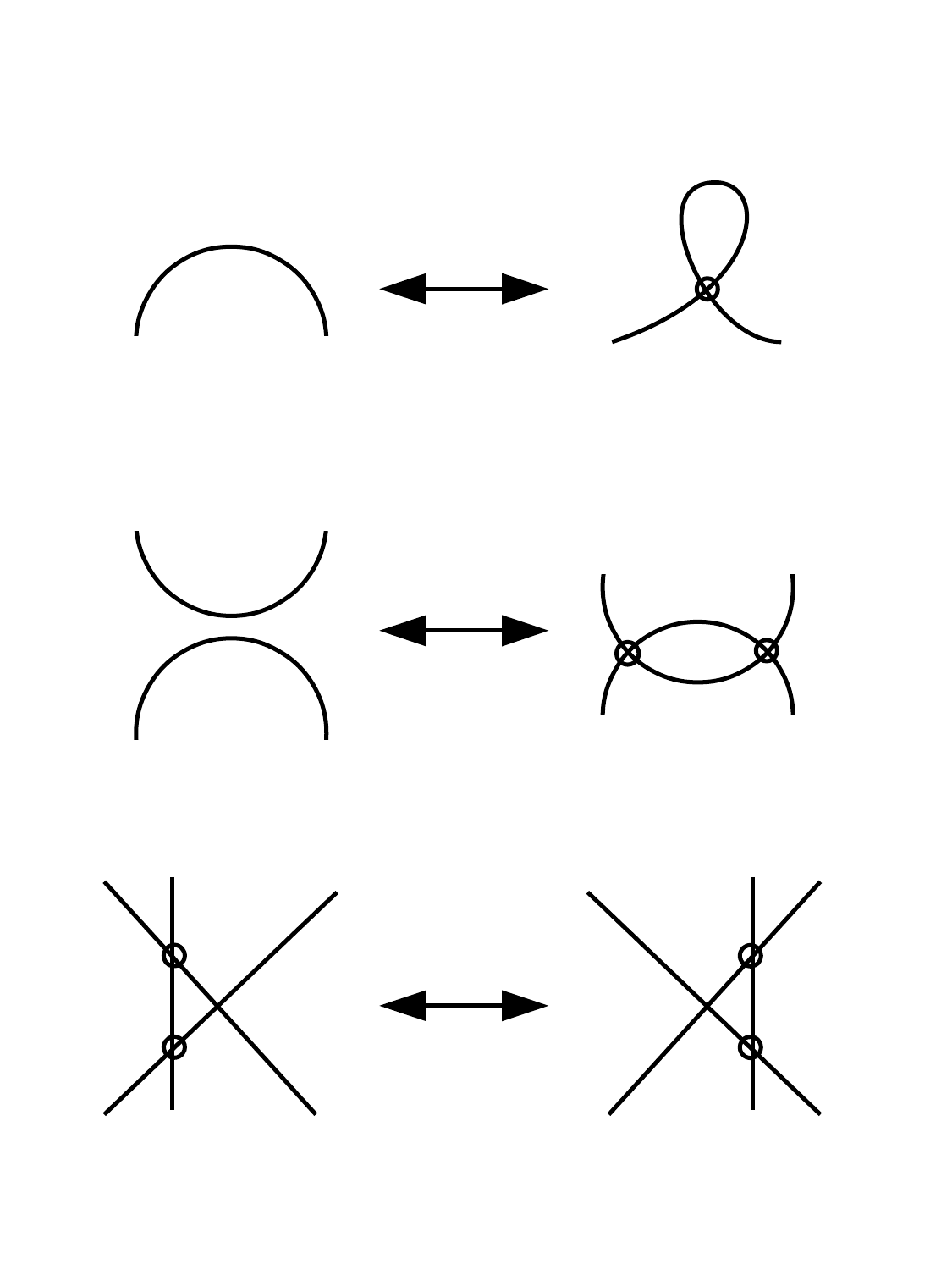}
		\caption{The additional virtual Reidemeister moves. Virtual crossings are indicated by crossings with small circles placed on them. The unmarked crossing in the third move can be either a virtual crossing or an arbitrary classical crossing.}
		\label{vReidemeister}
\end{figure}

Such Gauss codes may also be represented using \emph{virtual link diagrams}. A virtual link diagram is a link diagram where we allow arcs to cross one another in \emph{virtual} crossings, indicated in the diagram by a crossing marked with a little
circle. All the Reidemeister moves are permitted on such diagrams, and any arc segment with only virtual crossings may be replaced by a different arc segment with only virtual crossings, provided the endpoints remain the same. This is equivalent to allowing the additional moves shown in Fig. \ref{vReidemeister}.

Virtual links share many similarities with classical ones. For example, one can define virtual link groups and peripheral subgroups of virtual knot groups.
Furthermore, polynomial invariants of links, like the Alexander, Jones, and Homfly-pt polynomials, extend to virtual links.

\begin{figure}
		\centering
			\includegraphics[scale=0.5]{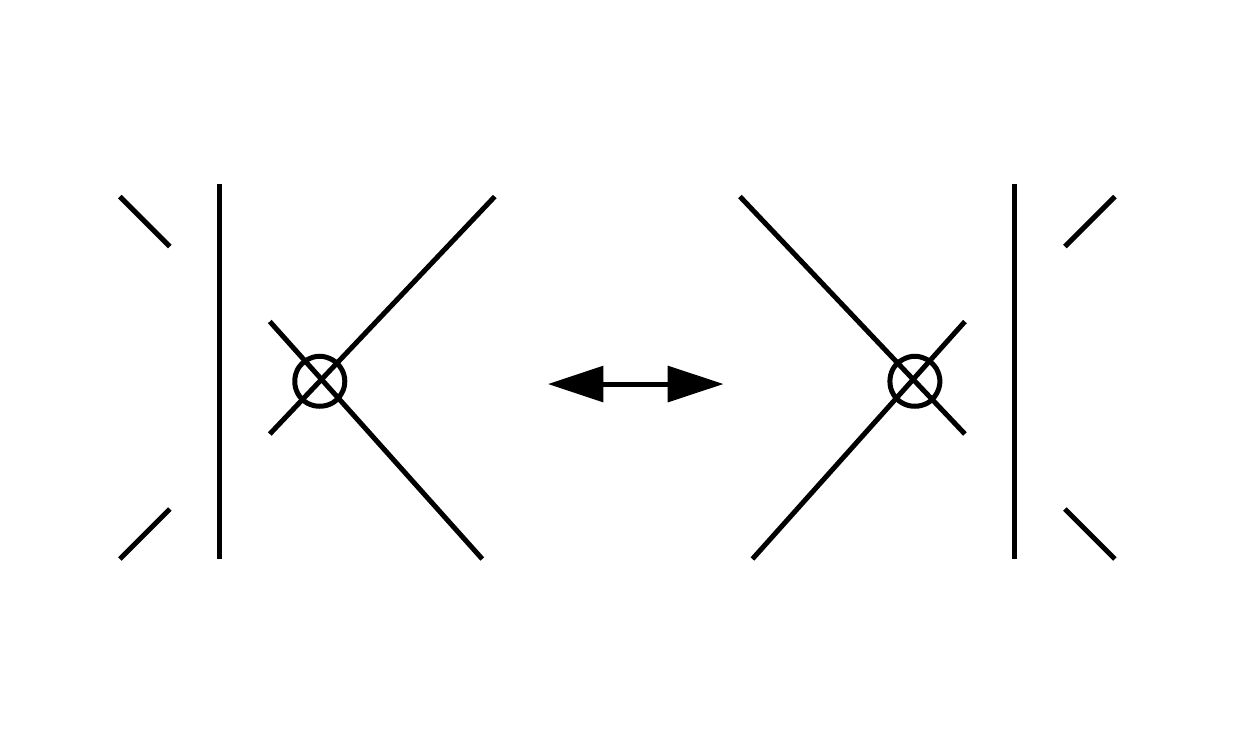}
		\caption{The ``forbidden move,'' which is not allowed for virtual $1$-links but is permitted for welded $1$-links.}
		\label{forbidden}
\end{figure}

The combinatorial definition and the geometric definition turn out to be equivalent to one another. There is a weaker equivalence relation on virtual link diagrams, which Satoh referred to as $w$-equivalence, but which is usually called \emph{welded} equivalence in the literature, from its original use in the welded permutation group (or welded braid group), \cite{FRR}. This is obtained from virtual equivalence on link diagrams by added in the \emph{welded forbidden move}, Fig. \ref{forbidden}. Virtual 1-links considered up to welded equivalence are usually just called welded links, a name which comes from its original use in \cite{FRR} for the welded braid group.

\section{Virtual $m$-links}

The definition of virtual $m$-links given in \cite{BKW} is a straightforward generalization of the geometric version of virtual 1-links given above. Let $K$ be an $m$-manifold embedded in $F\times I$, where $F$ is an $(m+1)$-manifold (possibly with boundary). Then $K$ is a virtual $m$-link. Just as before, define a relation $\sim$ by the condition that a link $K_1$ in $F_1 \times I$ is $\sim$-related to $K_2$ in $F_2\times I$ iff there exists an embedding $f:F_1\rightarrow F_2$ such that $f\times id_{[0,1]}(K_1)=K_2$. We define virtual equivalence to be the equivalence relation generated by $\sim$ together with smooth isotopy of links in thickened $(m+1)$-manifolds.

The knot group and knot quandle can be defined for virtual $n$-links, utilizing the \emph{Dehn space} of the virtual link. For a pair $(F\times [0,1], K)$, let $h(F)$ denote the space $(F\times [0,1])/r$, where $r$ is the relation  $(x, 1) r (x', 1)$ for all $x, x'\in F.$
Note that this is equivalent to gluing a cone over $F$ to $F\times\{ 1 \}$. Then the Dehn space of $K$, denoted $D(K)$, is given by $h(F)-K$. The importance of the Dehn space is the following theorem, which is proved in \cite{BKW}:

\begin{theorem}
The homotopy type of $D(K)$ is invariant under virtual equivalence.\label{hinv}
\end{theorem}

\begin{remark}
In fact, the homotopy type of $D(K)$ is invariant under welded equivalence of virtual 1-links as well. See \cite{BKW2}.
\end{remark}
We can define the link group $G(K)$ to be the fundamental group of $D(K)$. Similarly, we can define the link quandle $Q(K)$ to be the quandle of $D(K)$. See \cite{Mat, DJ} for the original definitions of quandles. In particular, the geometric definition of Joyce, \cite{DJ}, carries over to $D(K)$ with obvious modifications. Joyce's results showing that the quandle is equivalent to the fundamental group together with the peripheral subgroup and meridian also generalize to $D(K)$, \cite{BKW}.

\section{Generalizing virtual arcs}
\begin{figure}
		\centering
			\includegraphics[scale=0.5]{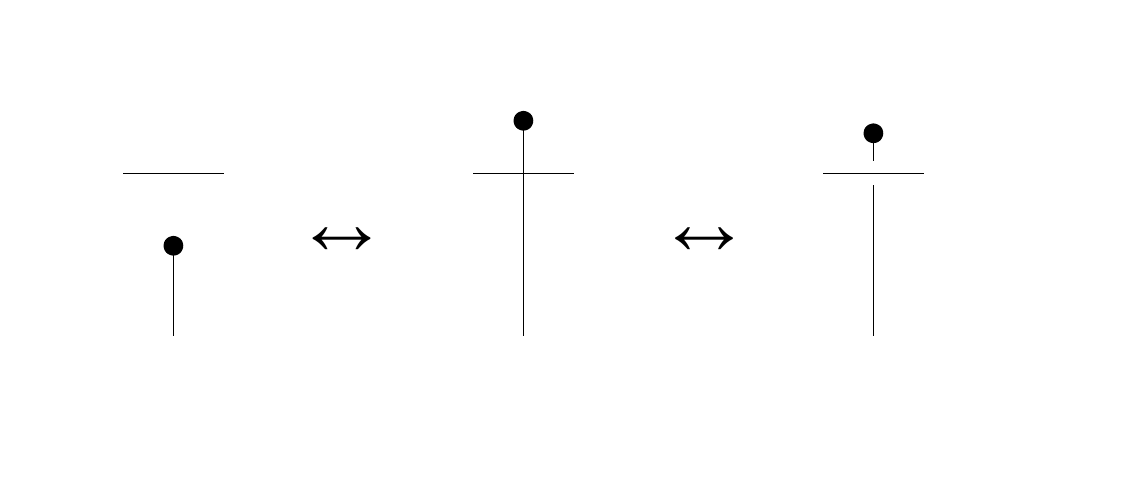}
		\caption{Allowed moves for ends of virtual arcs.}
		\label{ends}
\end{figure}

Satoh's virtual arcs are similar to the usual virtual links, but are allowed to have components with ends (necessarily two ends, since each component is the image of a 1-manifold). Endpoints are allowed to move under other arcs, but not over them, as in Fig. \ref{ends}.

Otherwise, link diagrams allowing arcs are the same as standard virtual links.

We can give a geometric definition of virtual arc diagrams. We will in fact give the definition that generalizes to higher dimensions. Because the arc diagrams bear some resemblance to "knotoids," we propose the name $m$-linkoid for these objects.

\begin{definition}
Let $F$ be an $(m+1)$ manifold, possibly with boundary. Then a virtual $m$-linkoid is an embedded $m$-manifold $K$ in $F \times I$, where $K$ may have a non-empty boundary, and we required that $\partial K$ be contained in the interior of $F \times \{ 0\}$. We define two virtual $m$-linkoids to be virtually isotopic if they are related by the following equivalence relation: Let $\sim$ be the relation defined by the condition that a link $K_1$ in $F_1 \times I$ is $\sim$-related to $K_2$ in $F_2\times I$ iff there exists an embedding $f:F_1\rightarrow F_2$ such that $f\times id_{[0,1]}(K_1)=K_2$. Then virtual equivalence, or isotopy, is the equivalence relation generated by $\sim$ together with isotopy of $K$ such that throughout the isotopy, $\partial K$ remains in the interior of $F \times \{ 0\}$.
\end{definition}

\begin{theorem}
The usual equivalence of geometric virtual 1-links with virtual 1-link diagrams extends to an equivalence between virtual 1-linkoids and virtual link diagrams with arcs.
\end{theorem}
The proof of this theorem is entirely straightforward.

\section{Satoh's Tube map}

\begin{figure}
		\centering
			\includegraphics[scale=0.5]{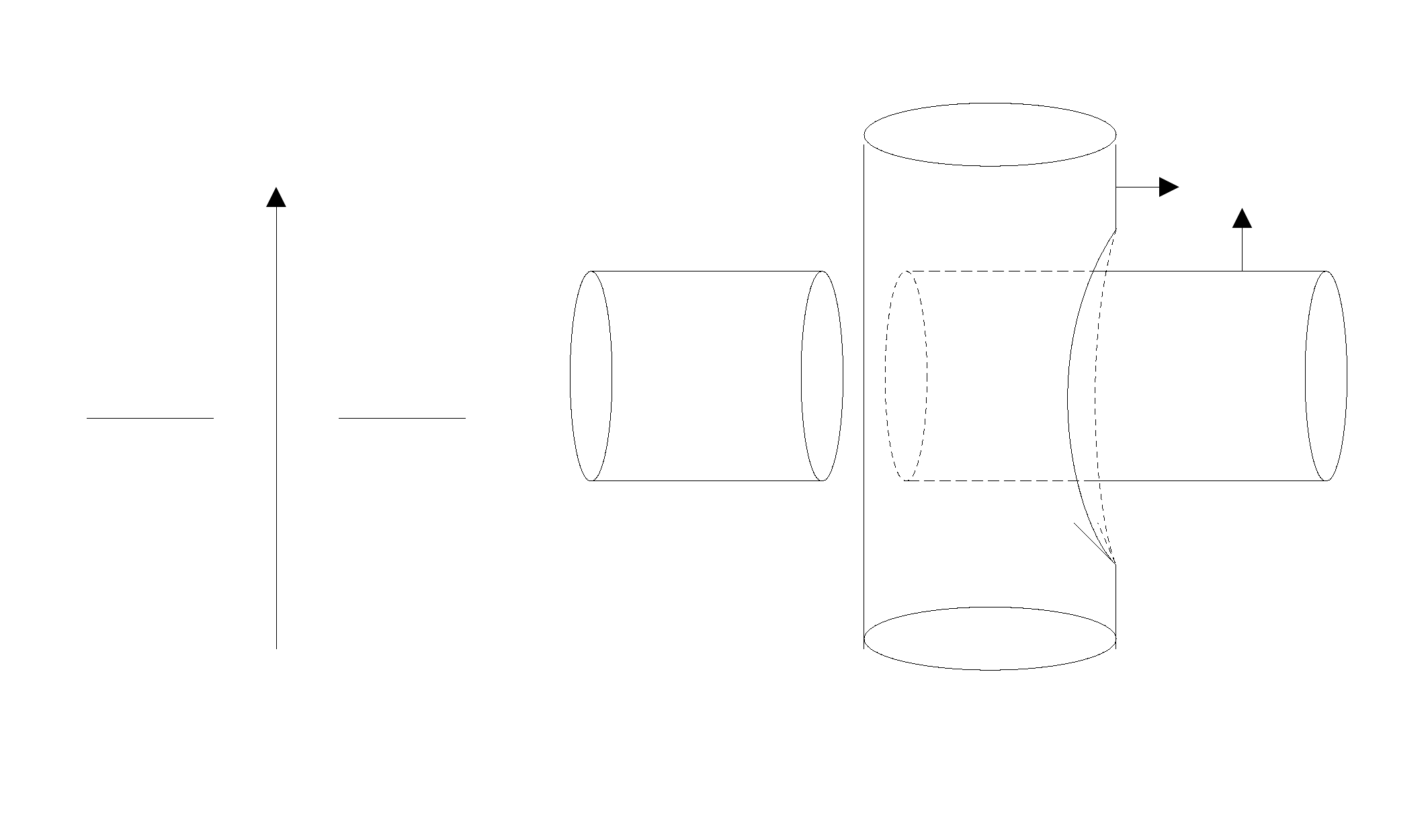}
		\caption{A piece of the broken surface diagram, right, corresponding to a crossing of the virtual link diagram, left.}
		\label{Tube}
\end{figure}

There is a map $Tube$ which takes virtual link diagrams to ribbon tori in $S^4$. The correspondence is highly geometric. Each arc becomes a cylinder. If an arc terminates (not terminates in a crossing, but just terminates), we cap the cylinder. At each crossing of the virtual link diagram, we connect the arcs as shown in Fig. \ref{Tube}. For a virtual link diagram $K$ we denote this ribbon torus as $Tube(K)$.

\begin{theorem}
If $K$ and $K'$ are $w$-equivalent, then $Tube(K)$ and $Tube(K')$ are isotopic through a stable equivalence of ribbon links.
\end{theorem}
This is shown in Satoh \cite{SS}. It is not hard to check each move.

\begin{theorem}
The link group and quandle of $Tube(K)$ is isomorphic to the link group and quandle of $K$.
\end{theorem}
This is also shown in \cite{SS}. The proof is straightforward: each cylinder in $Tube(K)$ gives rise to a generator corresponding to the generator from an arc in $K$. The relations at the crossings are the same, and the theorem follows.

\section{Generalizing to higher dimensions}
Our first observation is that we can generalize the notion of the $Tube$ map to higher dimensions. First, some notation. We will use $T_{mn}$ to denote the generalized $Tube$ map from $m$-links into $S^{n+2}$. The construction which we will use is not terribly complicated, but it is somewhat technical to write down explicitly.
Let $K:M\rightarrow F^{m+1}\times I$ be a virtual $m$-link, where $M$ has no boundary. Choose $n\geq 2m$. We will first build a $(n+1)$ CW complex $C$ from $K$. Start with $C_0 = M \times D^{n-m+1}$. Consider the projection of $K$ onto $F$. Then there will be a set of crossings in the projection, which lifts to a subspace of $K$. We create $C$ by applying the following equivalence relation: If $x_0, x_1 \in M$ project to the same point on $F$, and $K(x_1)$ is higher than $K(x_0)$ when projected onto the $I$ component, then we identify $\{ x_0\} \times D^{n-m+1}$ with a disk in the interior of $\{ x_1\} \times D^{m-n+1}$.

\begin{theorem}
$C$ embeds uniquely up to isotopy into $S^{n+2}$ as long as $n\geq 2m$.
\end{theorem}
\bpr
$C$ deformation retracts to a finite $m$-complex. All finite $m$-complexes embed uniquely up to isotopy in $S^{2m+2}$. But $C$ then embeds into a tubular neigborhood of this embedding, and any isotopy extends to this tubular neigborhood.
\epr

\begin{definition}
Given a virtual $m$-link $K$ and a choice of $n\geq 2m$, let $C$ be given as in the previous discussion. Let $i:C\rightarrow S^{n+2}$ be the embedding map. Finally, let $k:K\times S^{n-m}\rightarrow C$ be the restriction of the quotient map from $C_0$ to $C$ restricted to the boundary of $C_0$. We define $T_{mn}(K)$ to be the $n$-knot given by $i\circ k:K\times S^{n-m} \rightarrow S^{n+2}$.
\end{definition}

By the previous theorem, this construction is uniquely well-defined up to isotopy.

\begin{question}
This construction need not be well defined for $m<n<2m$. For a given $K$, can we characterize the minimal $n$ for which it is well-defined? In general this will depend upon $K$. For classical $K$, any $n\geq m$ will work, for example.
\end{question}

Intuitively, we are mapping a disk bundle over $M$ into $S^{n+2}$, subject to the following construction. The boundary must be embedded. Each disk fiber is either embedded or, if it meets a second disk fiber, then the image of one fiber is contained in the interior of the image of the other. If the bundle over $x_0$ is contained in the interior of the bundle over $x_1$, this corresponds to $K(x_0)$ being directly under $K(x_1)$ in the virtual link. Thus, crossings of the virtual $m$-link turn into fiber containment relations. It is instructive to consider this for Satoh's original map. A crossing from a virtual 1-link turns into one cylinder tube passing through the middle of another.

\begin{theorem}
$T_{mn}(K)$ is invariant under virtual isotopies of $K$. Indeed, if $K, K'$ are virtually isotopic, $T_{mn}(K)$ and $T_{mn}(K')$ are isoptic through an isotopy of knots all in the image of the $T_{mn}$ map, $n\geq 2m$.
\end{theorem}
\bpr
Parts of the isotopy that do not change the double point set of the projection of $K$ obviously satisfy this condition. Those which change the double point set change the complex $C$ in ways that are easily checked to change $T_{mn}(K)$ by an isotopy which leaves the special structure of bounding a $CW$ complex embedding unchanged, although the complex itself may change. These changes to the complex merely result in sliding $(n-m+1)$ disks through one another in $S^{n+2}$.
\epr

\begin{theorem}
The link group and quandle are preserved by $T_{mn}$.
\end{theorem}
\bpr
This follows from a straightforward application of the Van Kampen theorem. Each face of $K$ corresponds to a face of $T_{mn}(K)$, and the crossings correspond to a sphere of $T_{mn}(K)$ passing through another sphere. These give the same Wirtinger relations as those of $K$.
\epr

Furthermore, since $T_{mn}$ is invariant under virtual equivalence, it follows that we may use any invariants of classical $n$-links to distinguish virtual $m$-links.

\begin{question}
The fundamental group of the Dehn space of $K$ is isomorphic to the link group of $T_{mn}(K)$, which is the fundamental group of its complement. What is the relationship between the higher homotopy groups of the Dehn space of $K$, and the homotopy groups of complement of $T_{mn}(K)$?
\end{question}

\begin{remark}
We can extend all of the above to knotoids. If the underlying manifold $M$ of the virtual knot $K$ has a boundary, then we require that $\partial M \times D^{n-m+1}$ should not contain the images of any other disk fibers of $M\times D^{n-m+1}$, although their disk fibers may be contained in the interiors of the other disk fibers. The rest of the construction is similar. The boundary of $M$ corresponds to capping off the sphere fibers of $T_{mn}(K)$.
\end{remark}

We can also extend a result of Satoh about Artin's spinning construction. Note that for simplicity, we will discuss this construction in Euclidean space rather than in spheres. For a classical 1-link $K$, let $Spun(K)$ denote the torus in $\R^4$ produced in the following manner. Place $K$ in the upper half space $\R^3 _+$. Then form the link (whose components will be tori) in $\R^4$ by taking $(\R^3 _+, K)\times S^1$ quotiented by the relation that $(x, y, 0, t) ~ (x, y, 0, t')$ for any $t, t' \in S^1$. Satoh showed that $Spun(K)=Tube(K)$. Our version is somewhat more complicated.

We define $Sp_n (K)$ for a classical $m$-link $K$ in the following manner. Place $K$ in the upper half space $\R^{m+2} _+$, and form a link in $\R^{m+3}$ by taking $(\R^{m+2} _+, K)\times S^1$ quotiented by $(x_1, ... x_{m+1}, 0, t) ~ (x_1, ... x_{m+1}, 0, t')$ for any $t, t' \in S^1$. Notice that this $(m+1)$ link naturally bounds a CW complex of exactly the kind that we constructed above: in $(\R^{m+2} _+, K)$, we can consider extending a straight line down from each point of $K$. This line spins to form a 2-disk, and together, these disks form a CW complex. Then consider this as a sort of equatorial slice of an $n$-link in $\R^{n+2}$ bounding a similar CW complex, where the disks are thickened into $(n-m+1)$ disks. The resulting complex is identical, as a CW complex, to that formed when constructing $T_{mn}(K)$. But the choice of embedding of the CW-complex is irrelevant, as shown earlier. Hence:

\begin{theorem}
$Sp_n(K)\cong T_{mn}(K)$ for classical $K$.
\end{theorem} 

In fact, this leads to another method for constructiong the Dehn space $T_{mn}(K)$ using a virtual link, which is useful for studying the relationships between the Dehn spaces of $K$ and $T_{mn}(K)$. Let $K$ be a virtual link in $F\times I$. Thus, the image of $K$ is a subspace of $F\times I$; we abuse notation by writing $(F\times I, K)$ for this pair. Considered $S^{n-m}$ as a pointed space with basepoint $s_0$ and antipodal point $s_0 '$. Now consider $(F\times S^{n-m}\times I/r, K\times S^{n-m})$, where $r$ is the following relation: $(x, s, 0)r(x, s', 0)$ for any $x, s, s'$. We can therefore think of $(F\times S^{n-m}\times I/r)$ as the product of $F$ with a $D^{n-m+1}$ disk. We must reparametrize this into $F\times D^{n-m}\times I$ to make this a virtual link. To do so, we make a neighborhood of $F\times \{ s_0\} \times \{ 1\}$ into the "top" of the $D^{n-m}$, and a neighborhood of $F\times \{ s_0 '\} \times \{ 1\}$ into the "bottom" of the disk. Then we naturally complete the foliation. Denote the resulting virtual link as $T'_{mn}(K)$.

\begin{theorem}
If $K, K'$ are equivalent as virtual $m$-links, then $T'_{mn}(K), T'_{mn}(K')$ are equivalent as virtual $n$ links, for $n\geq m+1$.
\end{theorem}

Notice that $(F\times I, K)$ fits into this new paired subspace in a natural way, and that this implies there is a natural relationship between the Dehn spaces. Also, notice:

\begin{theorem}
$T'_{mn}(K)$ is virtually equivalent to $T_{mn}(K)$, when $n\geq 2m$.
\end{theorem} 
\bpr
Consider extending a straight line down from each point of $K$ in $F\times I$. This spins to give a CW-complex $C$ identical to that used to construct $T_{mn}(K)$. In addition, the projection of $C$ is the same as the projection of $K\times S^{n-m}$. Now consider $C$ embedded in $D^{n+1}\times I$. It can be isotoped so that fibers are all parallel to the $I$ component, and positioned so that in its projection, the images of two fibers meet iff one fiber is contained in the other (note: this is possible because we are working in $D^{n+1}\times I$ with $n\geq 2m$ and so $C$ embeds in $D^{n+1}$; otherwise, we could not guarantee this). We can therefore map a neighborhood of the projection of $C\subset D^{n+1}\times I$ to a neighborhood of the projection of $C$ in $T'_{mn}(K)$, since these neighborhoods are homeomorphic. In particular, two fibers project to the same point only if one is contained in the other, so the neighborhoods must be homeomorphic. But this induces a virtual equivalence.
\epr

This need not be the case when $n$ is not greater than $2m$, since in that case, $T_{mn}$ may not be uniquely defined.

The construction of $T'_{mn}$ should be compared to the vertical spin operation defined in \cite{BKW}, which is essentially our operation here when $n=m+1$. There is another spinning operation for virtual $m$-links: given $(F\times I, K)$, we could form the space $(F\times S^{n-m}\times I, K\times S^{n-m})$ and consider this as a virtual $n$-link, since $F\times S^{n-m}\times I$ is naturally foliated. In general this construction may not give us the same virtual link as $T'_{mn}$. We will refer to this procedure as a \emph{horizontal $(n-m)$-spin} of $K$.

\section{Generalizing $w$-equivalence}

We now turn to the question of $w$-equivalence for higher dimensions. In fact, we will see that it is potentially natural to generalize welded knots in such a manner that not all welded $n$-knots are necessarily virtual knots; at the least, they cannot be made into virtual knots by arbitrarily small perturbations.

Let $n\geq 2m$, and consider $M_0=M^m\times D^{n-m+1}$. We will say that a \emph{welded $m$-link} is a map $K:M_0 \rightarrow S^{n+2}$ which restricts to an embedding on the boundary, and such that if $K(m_0, x) = K(m_1, y)$, then either $K(\{ m_0\} \times  D^{n-m+1})$ is contained in the interior of $K(\{ m_1\} \times  D^{n-m+1})$, or 
$K(\{ m_1\} \times  D^{n-m+1})$ is contained in the interior of $K(\{ m_0\} \times  D^{n-m+1})$, such that the orientations of the disk fibers agree.

Intuitively, we are considering maps of $M_0$ to $S^{n+2}$, where $M_0$ is a trivial disk bundle over $M$. We require that the boundary of $M_0$ be embedded, and that if the images of two fibers meet, then the image of one fiber is contained in the interior of the image of the other fiber, with identical orientations. The fiber over $x_0$ being contained in the fiber over $x_1$ corresponds to $x_0$ passing "under" $x_1$.

\begin{definition}
Two welded $m$-links are \emph{$w$-equivalent} if they are connected by a map $L:M_0\times I \rightarrow S^{n+2}$ such that for each $x\in I$, $L(M_0\times \{ x\})$ is a welded $m$-link.
\end{definition}

\begin{remark}
There is a natural forgetful map from welded $m$-links to $n$-links, which comes from just considering the map on the boundary as a standard link.
\end{remark}

\begin{remark}
For $m=1$ our definition coincides with the usual notion of a welded link. This is easily seen by noting that our CW complex can generically be encoded as a Gauss code, and our welded equivalence corresponds to the welded moves on Gauss codes. In fact, our construction is rather similar to Rourke's construction for framed welded links, \cite{Rourke}. We could generalize this notion of framed links by requiring using $\R^{m+1} \times \R^{n+1-m}$ as the target space, and requiring that the disk fibers map into fibers of the target space, while the tangent space of the boundary should remain transverse to the fibers of the target space. If we do this, then instead of simply a quandle, we can define a fundamental link rack (see the next definition).
\end{remark}

\begin{definition}
We define the link group and quandle of a welded $m$-link to be the link group and quandle of the welded $m$-link considered as a classical $n$-link.
\end{definition}

\begin{theorem}
The link group and quandle are invariants of welded $m$-links, and also do not depend on the choice of $n\geq 2m$.
\end{theorem}

Now, we must remark that welded $m$-links, for $m\geq 2$, do not appear to necessarily correspond to virtual $m$-links. This is because of the possibility of sliding undercrossings past one another with welded equivalence. Consider a welded 2-link which locally has a region which, for some choice of codimension-$1$ foliation, can be interpreted as an isotopy of a welded $1$-link undergoing a forbidden move. There is no way to get this specific embedding from the generalized Tube maps $T_{2, n}$ we defined above.

However, every virtual $m$-link can be interpreted as a welded $m$-link via the $T_{mn}(K)$ construction.

\begin{question}
Is every welded $m$-link isotopic to the image of some virtual $m$-link under $T_{mn}$? We conjecture that the answer is no.
\end{question}

Note that there are examples of welded $m$-links which do come from virtual $m$-links under the tube map which are isotopic, despite the virtual links not being virtually equivalent. For example, consider creating a virtual $m$-link by horizontal 1-spinning the virtual trefoil $m$ times. This is a nontrivial virtual $m$-link, as may be seen from its stack invariants, \cite{BKW2}. However, the welded isotopy which makes the virtual trefoil $w$-equivalent to the unknot also spins to give a welded isotopy of the $m$-spun virtual trefoil, turned into a welded $m$-knot via the $T_{mn}$ map, to an unknotted $m$-knot.

\begin{remark}
We can extend all of the above to knotoids just as we did with the $T_{mn}$ map. If the underlying manifold $M$ of the virtual knot $K$ has a boundary, then we require that $\partial M \times D^{n-m+1}$ should not contain the images of any other disk fibers of $M\times D^{n-m+1}$, although their disk fibers may be contained in the interiors of the other disk fibers. This gives rise to the notion of a welded $m$-knotoid.
\end{remark}

\section{Another interpretation for virtual $m$-links}

We may in fact strengthen the notion of welded $m$-links given above to give a definition which is equivalent to the virtual $m$-links of \cite{Kup, BKW}.

Let $n\geq 2m$, and consider $M_0=M^m\times D^{n-m+1}$. Now consider maps $K:M_0 \rightarrow S^{n+2}$ which restricts to an embedding on the boundary, and such that if $K(m_0, x) = K(m_1, y)$, then either $K(\{ m_0\} \times  D^{n-m+1})$ is contained in the interior of $K(\{ m_1\} \times  D^{n-m+1})$, or 
$K(\{ m_1\} \times  D^{n-m+1})$ is contained in the interior of $K(\{ m_0\} \times  D^{n-m+1})$, such that the orientations of the disk fibers agree, with one additional restriction. We will require that if one disk fiber $D_0$ maps into the image of another disk fiber $D_1$, then the preimage of their intersection should form a disk in $D_1$ with radius greater than half the radius of $D_1$. Such a welded knot will be said to be in \emph{virtual tube form}.

This restriction prevents things like the forbidden move, since two disk fibers $D_1, D_2$ cannot be inside the same fiber $D_0$ unless one of them is contained in the other. In fact, it is not hard to check that:

\begin{theorem}
Every welded knot in virtual tube form is in the image of $T_{mn}$. Furthermore, $K, K'$ being virtual $m$-links which are virtually equivalent, is equivalent to $T_{mn}(K)$ and $T_{mn}(K')$ being equivalent through welded $m$-links in virtual tube form.
\end{theorem}

The theorem follows by observing that the containment relations of the disk fibers corresponds to double point sets for a virtual $m$-link.

It follows that an alternative definition of virtual $m$-links, equivalent to that given in \cite{BKW}, is welded $m$-links in virtual tube form, up to welded equivalence through other welded links in virtual tube form.

It may be wondered whether this result, which is rather inelegant from a certain perspective, is really all that useful. However, putting geometric restrictions on topologically knotted objects is not a new idea. Legendrian knots, for example, are knots that obey certain geometric restrictions with respect to a contact structure.

One rather immediate consequence of this result is the following theorem:

\begin{theorem}
If $K, K'$ are virtual $m$-links which are $w$-equivalent, then $T_{mn}(K)$ is isotopic to $T_{mn}(K')$, $n\geq 2m$.
\end{theorem}
\bpr
Interpret $K, K'$ as welded $m$-links in virtual form. The $w$ equivalence between them induces an isotopy of the links as $n$-links.
\epr

\begin{question}
Using this new characterization of virtual $m$-links, is it possible to construct new invariants or give new definitions for known invariants? For example, biquandle invariants?
\end{question}

Since biquandle invariants are only understood combinatorially, it would be useful to have a geometric definition for them, which we might hope could be found by using this characterization.

\section{Yet another interpretation of virtual $m$-links}

Consider the vertical double $VD(K)$ of a virtual $m$-link $K$. We form this as follows: given $K:M \rightarrow F\times I$, a virtual $m$-link, we isotope the map $K$ so that the image lies in $F\times (0.5, 1)$. Let $\pi_F, \pi_I$ be the natural projection of $F\times I$ to $F$ and to $I$ respectively. Then we form a new $m$-link $VD(K): M\times \{ 0, 1\} \rightarrow F\times I$ as follows: $VD(K)(m, 0)=K(m)$, $VD(K)(m, 1) = (\pi_F (K(m)), 1-\pi_I (K(m)))$. The result is that we adjoin to $K$ a copy of the vertical mirror image $K^*$ of $K$ directly underneath it. For convenience, we may wish to reverse the orientation of this copy. The vertical double was introduced in \cite{BKW2}.

Now interpret $VD(K)$ as a welded $m$-link with extra structure: observe that the disk fibers (of the CW complex we use in the $T_{mn}$ construction) of $VD(K)$ are naturally paired, with every fiber of $K^*$ corresponding to $VD(K)(m, 1)$ contained in the fiber of $K$ coming from $VD(K)(m, 0)$. Furthermore, if the fiber from $VD(K)(m_0, 0)$ contains the fiber from $VD(K)(m_1, 0)$, then the fiber from $VD(K)(m_1, 1)$ contains the fiber from $VD(K)(m_0, 1)$. A welded link whose base manifold is $M\times \{ 0, 1\}$ will be called a welded link in \emph{stack form} if it retains these containment relations.

\begin{theorem}
$K \cong K'$ as virtual $m$-links iff $VD(K)$ and $VD(K')$ are equivalent as welded $m$-links via a welded equivalence through welded links in stack form.
\end{theorem}
\bpr
The only if direction is trivial. For the if direction, observe that the changes to the containment structure on $VD(K)$ restricted to $M\times \{ 0\}$ must correspond to changes in the double point set of $K$ which are permitted by virtual equivalence. This is because undercrossings can no longer pass by one another, as is easily checked. Therefore, this equivalence induces an equivalence of $K$ and $K'$ as virtual $m$-links, as interpreted in the previous section.
\epr

\begin{question}
Can this condition be weakened to simply saying that $VD(K)$ and $VD(K')$ are equivalent as welded $m$-links (possibly as $m$-links in stack form), without requiring the welded equivalence to retain the stack form throughout the equivalence?
\end{question}

\section{Extending extended welded links}

The notion of an extended welded 1-link was introduced to incorporate the possibility of \emph{wens} in ribbon 2-links. See for example \cite{Damiani}. A wen corresponds to turning the disk fibers in the $Tube$ map upside down.

It is worth remarking that in our above discussion of welded $m$-links, if we consider these presented via Fox-Milnor movies of $(m-1)$ links, that although it is possible for welded crossings to move through undercrossings in general, a welded cusp cannot move through undercrossings. This peculiarity is alleviated if we modify our definition slightly, in a way which naturally generalizes the notion of an extended welded link.

Let $n\geq 2m$, and consider $M_0=M^m\times D^{n-m+1}$. We will say that an \emph{extended welded $m$-link}, or ewe-link, is a map $K:M_0 \rightarrow S^{n+2}$ which restricts to an embedding on the boundary, with the following additional condition: if $K(m_0, x) = K(m_1, y)$, where $x\in \partial D^{n-m+1}$, then there is a disk neighborhood $N$ of $(m_1, y)$, such that $K(\{m_0 \times D^{n-m+1}) \subset K(N)$. Two ewe $n$-links are \emph{ewe-equivalent} if they are connected by a map $L:M_0\times I \rightarrow S^{n+2}$ such that for each $x\in I$, $L(M_0\times \{ x\})$ is a ewe-$m$-link.

\begin{figure}
		\centering
			\includegraphics[scale=0.5]{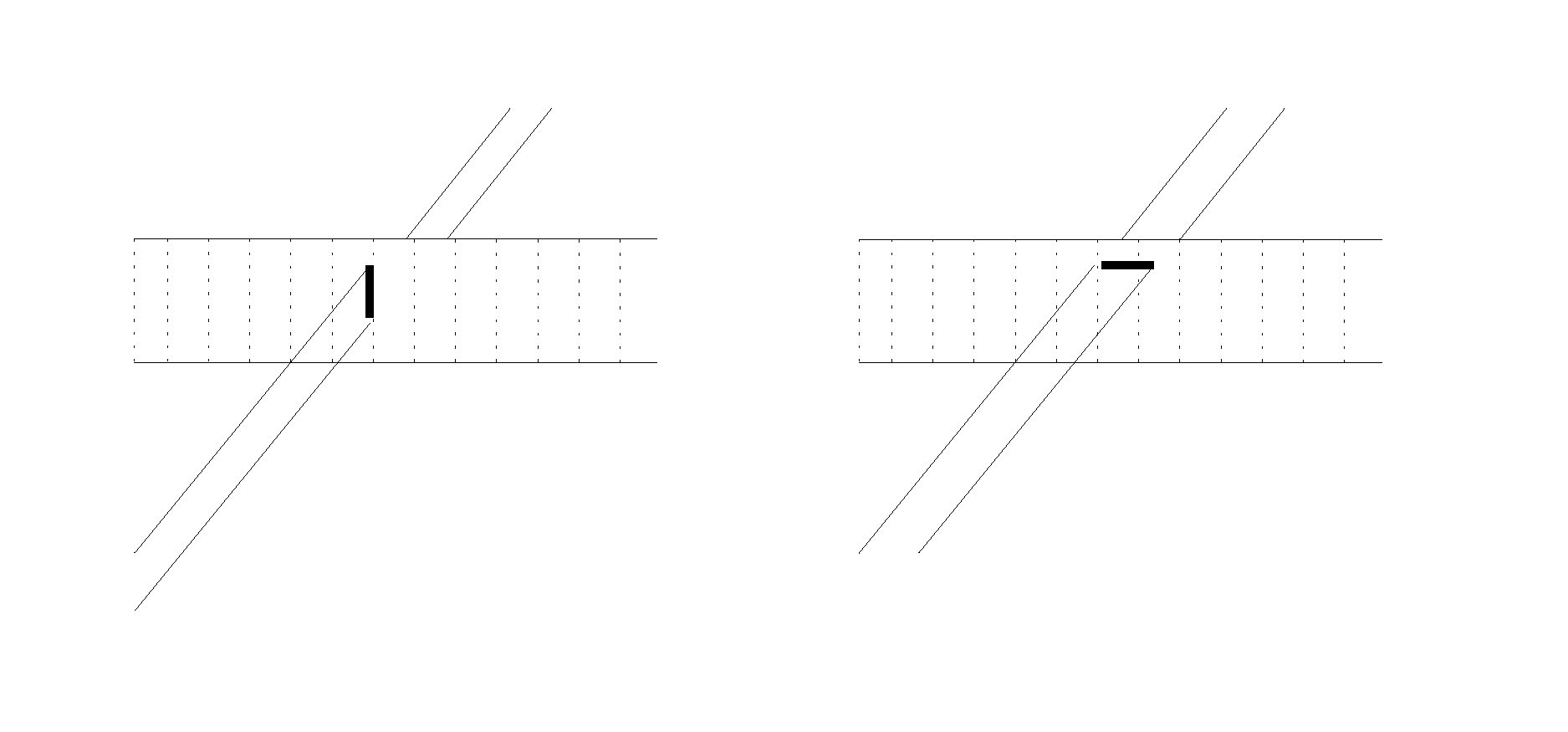}
		\caption{In this figure, we show a 3-dimensional slice of a welded $m$-link, left, and a ewe $m$-link, right. The situation on the right is not allowed for regular welded $m$-links, because the "undercrossing" fiber is not contained in any single "overcrossing" fiber. It is, however, contained in a disk neighborhood of an "overcrossing" fiber.}
		\label{ewe}
\end{figure}

What this allows is for the "undercrossing" disk fiber to not be embedded inside a single "overcrossing" disk fiber, but to be turned transversely to it. We illustrate this in Fig. \ref{ewe}.

The following theorem is rather immediate:

\begin{theorem}
Ewe 1-links are equivalent to extended welded 1-links as defined in \cite{Damiani}.
\end{theorem}

Thus, ewe $m$-links naturally generalize extended welded links to higher dimensions. And as before, we can naturally extend these notions to ewe $m$-knotoids.


\end{document}